\documentclass[12pt]{article}
\usepackage{latexsym,color,amsmath,amsthm,amssymb,amscd,amsfonts}

\newtheorem{lemma}{Lemma}[section]
\newtheorem{proposition}[lemma]{Proposition}
\newtheorem{remark}[lemma]{Remark}
\newtheorem{theorem}[lemma]{Theorem}
\newtheorem{definition}[lemma]{Definition}
\newtheorem{corollary}[lemma]{Corollary}

\newtheorem{conj}[lemma]{Conjecture}

\newtheorem*{remark*}{Remark}

\def\R{{\mathbb R}}

\def\eps{\varepsilon}

\makeatletter \@addtoreset {equation}{section}
\renewcommand\theequation
  {\ifnum \c@subsection>\z@ \arabic{section}.\arabic{subsection}.\arabic{equation}
  \else \arabic{section}.\arabic{equation} \fi}
\makeatother

\begin{document}

\title {On Symplectic Capacities and Volume Radius}
\author{Shiri Artstein-Avidan
\thanks{
The first named author was supported by the National Science
Foundation under agreement No. DMS-0111298. Any opinions, findings
and conclusions or recommendations expressed in this material are
those of the authors and do not necessarily reflect the views of the
National Science Foundation.} \, and Yaron Ostrover}
\maketitle

\noindent{\em {\bf Abstract:} In this work we present an improvement
to a theorem by C. Viterbo, relating the symplectic capacity of a
convex body and its volume. This provides one more step towards the
proof of the following conjecture: among all convex bodies in
${\R}^{2n}$ with a given volume, the Euclidean ball has maximal
symplectic capacity. More precisely, the conjecture states that the
best possible constant $\gamma_n$ such that for any choice of a
symplectic capacity $c$ and any convex body $K \subset {\R}^{2n}$ we
have
$$ {\frac {c(K)} {c(B^{2n})}} \leq \gamma_n \biggl ({\frac {{\rm
Vol}(K)} {{\rm Vol}(B^{2n})}} \biggr )^{1/n}$$ is $1$. Until this
work, the best estimate known to hold for general convex bodies,
coming from Viterbo's work was $\gamma_n={32n}$, and $\gamma_n={2n}$
in the case of centrally symmetric bodies. Our main result in this
text is that there exists a universal constant $A$ for which
$\gamma_n \le A \log^2(n)$ for all convex bodies in $\R^{2n}$.
Moreover, we show wide classes of convex bodies for which the
inequality holds without the logarithmic term. }

\section{Introduction and Main results}

This paper lies at the meeting point of Asymptotic Geometric
Analysis and Symplectic Geometry. In particular we use methods from
Asymptotic Convex Geometry (sometimes called the Local Theory of
Banach Spaces), to improve a result of Viterbo concerning symplectic
capacities of convex bodies. These methods are linear in nature, and
the reader should not expect any difficult symplectic analysis.
However, we stress that the naive linear approach provides a
significant improvement to the known results. The type of
improvement we provide is a reduction from order $n$ to order
$\log(n)$ of a dimension-dependent isoperimetric constant, where $n$
is the dimension of the space involved. Clearly, this reduction
becomes especially relevant in large dimensions. Understanding the
behavior of convex bodies in dimension tending to infinity is the
main subject of Asymptotic Geometric Analysis, which we here join
with the symplectic notion of capacity.

Consider the $2n$-dimensional Euclidean space ${\mathbb R}^{2n}$
with the standard linear coordinates $(x_1,y_1, \ldots,x_n,y_n)$.
One equippes this space with the standard symplectic structure
$\omega_{st} = \sum_{j=1}^n dx_j \wedge dy_j$, and with the standard
inner product $g_{st} = \langle \cdot,\cdot \rangle$. Note that
under the identification (see notations below) between ${\mathbb
R}^{2n}$ with $ {\mathbb C}^n$ these two structures are the real and
the imaginary part of the standard Hermitian inner product in
${\mathbb C}^n$. In this work we consider the class of convex bodies
in ${\mathbb R}^{2n}$. We are interested in comparing the symplectic
way of measuring the size of a convex body, using what is called
``symplectic capacities", with the standard Riemannian way limited
here to volume. In order to make this more precise we need some
preliminaries.

\begin{definition} \label{Def-sym-cap}
A 
symplectic capacity on $({\mathbb
R}^{2n},\omega_{st})$ associates to each  subset $U \subset
{\mathbb R}^{2n}$ a non-negative number $c(U)$ such that the
following three properties hold:
\begin{enumerate}
\item[(P1)] $c(U) \leq c(V)$ for $U \subseteq V$ (monotonicity)
\item[(P2)] $c \big (\psi(U) \big )= |\alpha| \, c(U)$ for  $\psi
\in {\rm Diff} ( {\mathbb R}^{2n} )$ such that $\psi^*\omega_{st}
= \alpha \, \omega_{st}$ (conformality)
\item[(P3)] $c \big (B^{2n}(r) \big ) = c \big (B^2(r) \times {\mathbb
C}^{n-1} \big ) = \pi r^2$ (nontriviality and normalization),
\end{enumerate}
\end{definition}
\noindent where $B^{2k}(r)$ is the open $2k$-dimensional ball of
radius $r$. Note that the third property disqualifies any
volume-related invariant, while the first two properties imply
that every two sets $U,V \subset {\mathbb R}^{2n}$ will have the
same capacity provided that there exists a symplectomorphism
sending $U$ onto $V$. Recall that a {\it symplectomorphism} of
${\mathbb R}^{2n}$ is a diffeomorphism which preserves the
symplectic structure i.e., $\psi \in {\rm Diff} ( {\mathbb R}^{2n}
)$ such that $\psi^* \omega_{st} = \omega_{st}$. We will denote by
${\rm Symp}({\mathbb R}^{2n}) = {\rm Symp}({\mathbb
R}^{2n},\omega_{st})$ the group of all the symplectomorphisms of
$({\mathbb R}^{2n},\omega_{st})$.

A priori, it is not clear that symplectic capacities exist. The
celebrated non-squeezing theorem of Gromov~\cite{G} shows that for
$R
> r$ the ball $B^{2n}(R)$ does not admit a symplectic embedding
into the symplectic cylinder $Z^{2n}(r):= B^2(r) \times {\mathbb
C}^{n-1}$. This theorem led to the following definitions:

\begin{definition} The symplectic radius of a non-empty set $U \subset {\mathbb R}^{2n}$ is
$$ c_B(U) := \sup \left \{\pi r^2 \, | \,
\ There \ exists \  \psi \in {\rm Symp}({\mathbb R}^{2n}) \ with \
\psi \left (B^{2n}(r) \right ) \subset U \right  \}.$$ The
cylindrical capacity of $U$ is
$$ {c}^Z(U) := \inf \left \{\pi r^2 \, | \,
\ There \ exists \  \psi \in {\rm Symp}({\mathbb R}^{2n}) \ with \
\psi (U) \subset Z^{2n}(r)  \right  \}.$$
\end{definition}

Note that both the symplectic radius and the cylindrical capacity
satisfy the axioms of Definition~\ref{Def-sym-cap} by the
non-squeezing theorem. Moreover, it follows from
Definition~\ref{Def-sym-cap} that for every symplectic capacity
$c$ and every open set $U \subset {\mathbb R}^{2n}$ we have
$c_B(U) \le c(U) \le c^Z(U)$.

The above axiomatic definition of symplectic capacities is
originally due to Ekeland and Hofer~\cite{EH}. Nowadays, a variety
of symplectic capacities can be constructed in different ways. For
several of the detailed discussions on symplectic capacities we
refer the reader
to~\cite{CHLS},~\cite{Ho},~\cite{HZ},~\cite{L},~\cite{Mc}
and~\cite{V1}.

In this work we are interested in an inequality relating the
symplectic capacity of a convex body in ${\R}^{2n}$ and its
volume. By a convex body we mean a convex bounded set in $\R^{2n}$
with non-empty interior. This inequality supports the conjecture
that among all convex bodies in ${\R}^{2n}$ with a given volume,
the symplectic capacity is maximal for the Euclidean ball. Note
that by monotonicity this is obviously true for the symplectic
radius $c_B$. More precisely, denote by ${\rm Vol}(K)$ the volume
of $K$ and abbreviate $B^{2n}$ for the open Euclidean unit ball in
${\mathbb R}^{2n}$. Following Viterbo~\cite{V} we state

\begin{conj} \label{The-conjecture} For any symplectic capacity $c$ and for any convex
body $K \subset {\mathbb R}^{2n}$
$${\frac {c(K)} {c(B^{2n})}} \leq
 \biggl ( {\frac {{\rm Vol}(K)} {{\rm Vol}(B^{2n})} } \biggr
)^{1/n}$$ and equality is achieved only for symplectic images of
the Euclidean ball.
\end{conj}

The first result 
in this direction is due to Viterbo~\cite{V}. Using John's ellipsoid
(which also evolved in Convex Geometric Analysis) he proved:
\begin{theorem}[Viterbo]\label{cnjV}
For a convex body $K \subset {\mathbb R}^{2n}$ and a symplectic
capacity $c$ 
$$ {\frac {c(K)} {c(B^{2n})}} \leq
\gamma_n   \biggl ( {\frac {{\rm Vol}(K)} {{\rm Vol}(B^{2n})} }
\biggr )^{1/ n}$$ where $\gamma_n = {2n}$ if $K$ is centrally
symmetric and $\gamma_n = 32n$ for general convex bodies.
\end{theorem}

In~\cite{H}, Hermann constructed a starshaped domain in ${\mathbb
R}^{2n}$, for $n > 1$, with arbitrarily small volume and some fixed
cylindrical capacity. Therefore, in the category of starshaped
domains the above theorem with {\em any} constant $\gamma_n$
independent of the body $K$, must fail. In addition, he proved the
above conjecture for a special class of convex bodies which admit
many symmetries, called convex Reinhardt domains (for definitions
see~\cite{H}).

Here we provide one more step towards the proof of the above
conjecture. Before we state our main results we wish to emphasize
that we work exclusively in the category of linear symplectic
geometry. That is, we restrict ourselves to the concrete class of
linear symplectic transformations. It turns out that even in this
limited category, there are tools which are powerful enough to
obtain a significant improvement of Theorem~\ref{cnjV} above.
More precisely, let ${\rm Sp}({\mathbb R}^{2n})={\rm Sp}({\mathbb
R}^{2n},\omega_{st})$ denote the group of linear symplectic
transformation of ${\mathbb R}^{2n}$. We consider a more
restricted notion of linearized cylindrical capacity, which is
similar to $c^Z$ but where the transformation $\psi$ is taken only
in ${\rm Sp}({\mathbb R}^{2n})$ namely
$$ {c}_{lin}^Z(U) := \inf \left \{\pi r^2 \, | \,
\ There \ exists \  \psi \in {\rm Sp}({\mathbb R}^{2n}) \ with \
\psi (U) \subset Z^{2n}(r)  \right  \}.$$ Of course, it is always
true that for every symplectic capacity $c$ we have $c \le c^Z \le
c_{lin}^Z$.

Our main result is the following

\begin{theorem}\label{The_Main_Theorem}
There exists a universal constant $A_1$ such that for every even
dimension $2n$ and any convex body $K \subset {\mathbb R}^{2n}$ 
we have
$$ 
{\frac {c^Z_{lin}(K)} {c(B^{2n})}}  
\leq A_1 \log^2(n)
 \biggl (
{\frac {{\rm Vol} (K)} {{\rm Vol} (B^{2n})} } \biggr )^{ 1/ n}.$$
\end{theorem}
\noindent The theorem clearly implies that Theorem \ref{cnjV} holds with
$\gamma_n = A_1 \log^2(n)$. We remark that the methods Viterbo
used to prove Theorem \ref{cnjV} were also linear.

In addition we show that there are certain, quite general, families
of convex bodies for which the above inequality is true without the
logarithmic factor. Such examples are the unit balls of $\ell_p^n$
for $1 \le p \le \infty$, all {\it zonoides} (symmetric convex
bodies which can be approximated by Minkowski sums of line segments
in the Hausdorff sense, in particular all the projections of cubes),
and all bodies which satisfy some more complicated geometric
conditions called bounded type-$p$ constant for $p>1$. The precise
definitions and results for special classes of convex bodies will be
given in Section~\ref{IMP} below.

Some of the above mentioned examples are based on a strengthened
formulation of Theorem~\ref{The_Main_Theorem}. This strengthening is
attained by bounding from above the symplectic capacity of a convex
body $K$ by some parameter of the body which measures, roughly
speaking, the ``difference" between the Banach space $X_K$ whose
unit ball is the body $K$ and a  Hilbert space. In order to be more
precise we need to introduce the notion of the $K$-convexity
constant of a Banach space $X$, known also as the Rademacher
projection constant. We begin with the following preliminaries.
Consider the group $\Omega^m =  ({\mathbb Z} / {2 { \mathbb Z}})^m
\simeq \{-1,1\}^m$ and let $\mu$ denote the uniform probability
measure on $\Omega^m$ (i.e., normalized counting measure). For each
$A \subseteq \{1,\ldots,m\}$ we define the Walsh function $W_A \in
L^2(\Omega^m,{\mathbb R})$ by
\[ W_A(t) = \prod_{i \in A}r_i(t) = \prod_{i \in A} t_i,\] where $r_i : \Omega^m
\rightarrow \{-1,1 \}$ are the Rademacher functions $\ r_i(t) =
t_i$, i.e., the $i$th coordinate function (so, $r_i = W_{\{i\}}$).
We set $W_{\emptyset} \equiv 1$. Note that $|W_A|=1$ for all $A
\subseteq \{1,\ldots,m\}$, and for every $A \neq B$ the functions
$W_A, W_B$ are orthogonal, namely
\[ \langle W_A,W_B \rangle = 
\int_{\Omega^m} \prod_{i \in A} r_i \prod_{j
\in B}r_j \, d\mu= \int_{\Omega^m} \prod_{k \in A \triangle B} r_k
\, d\mu = \prod_{k \in A \triangle B} \int_{\Omega^m} r_k \, d\mu
= 0
\]

The Walsh functions form an orthonormal basis of $L^2(\Omega^m,
\R)$. Alternatively, $(W_A)_{A \subseteq \{1, \ldots, m\}}$ is the
group of characters of the multiplicative group $({\mathbb Z} /
{2{\mathbb Z}})^m$. For a Banach space $X$ define
\[ L^2(\Omega^m,X) = \Bigl \{ f : \Omega^m \rightarrow X \ ; \
\|f\|_{L^2(\Omega^m,X)} = \Bigl ({\int_{\Omega^m}}\|f\|^2_X \, d\mu
\Bigr )^{1/2} \Bigr \}. \]
%
The space $L^2(\Omega^m,X)$ is always a Banach space, and it is a
Hilbert space if and only if $X$ is. Still, every $f \in
L^2(\Omega^m,X)$ can be represented as
\[ f = \sum_{A \subseteq \{1,\ldots,m\}} {\widehat f}(A) W_A, \  {\rm
where} \ \ {\widehat f}(A) = \int_{\Omega^m}fW_A \, d\mu = {\frac
1 {2^m}} \sum_{t \in \Omega^m}f(t)W_A(t).\]
We will be interested in one special subspace of $L^2(\Omega^m,X)$
namely the one spanned by the $m$-Rademacher functions: \[ Rad_mX
= \bigl \{ \sum_{i=1}^m x_i r_i \ ; \ x_i \in X \ , \ i=1,\ldots,m
\bigr \}, \] equipped with the $L^2$-norm. Consider the Rademacher
projection operator ${\cal R}_m : L^2(\Omega^m,X) \to Rad_mX $
defined by 
\[ {\cal R}_m (f) =
\sum_{i=1}^m {\widehat f}({\{i\}})r_{i}.
\]
\begin{definition} The Rademacher projection constant of a Banach space
$X$ is the supremum of the operator norms of the projections
${\cal R}_m$, i.e.
\[ \|{\rm Rad}\|_X := \sup_m \|{\cal R}_m \|.\]
\end{definition}

This is also known in the literature as the $K$-convexity constant
of $X$. Note that for an infinite dimensional Banach space $X$,
this constant may be infinite. When for some (infinite
dimensional) Banach space $X$ the number $\|{\rm Rad}\|_X = C$ is
finite,
$X$ is called $K$-convex with $K$-convexity constant $C$. In what
follows we shall deal only with finite dimensional spaces,
and it is not difficult to check that they are always $K$-convex.
However,
the dependence of the $K$-convexity constant on the dimension of
the space will be of interest to us.
For example, if we consider a family of finite dimensional Banach
spaces $X_k$
(with increasing dimension, say) which all
arise as subspaces of some fixed infinite dimensional Banach space $X$
which is $K$-convex, then we know that the $K$-convexity constants
of the spaces $X_k$ are uniformly bounded by the $K$-convexity constant
of $X$.

In the special case where $X$ is a Hilbert space, each ${\cal R}_m$
is an orthogonal projection of norm equal to $1$. A fundamental
theorem of Pisier, see \cite{Pi2}, states that for every Banach
space $X$ which is isomorphic to a Hilbert space $H$,
\[ \|{\rm Rad}\|_X \leq c \log [d(X,H) + 1],\]
where $c$ is a universal constant and $d$ denotes the Banach-Mazur
distance, defined for two isomorphic normed spaces as $\inf \{ \|T\|
\, \|T^{-1}\|\}$ where the infimum runs over all isomorphisms $T:X
\rightarrow H$. Notice that for a $k$-dimensional $X$, John's
Theorem (see \cite{MS}) implies that $d(X,\ell_2^k) \leq \sqrt{k}$.
Thus, combining John's Theorem and Pisier's result we obtain that
Theorem~\ref{The_Main_Theorem} actually follows from

\begin{theorem}\label{The_Main_Theorem-after-rad}
There exists a universal constant $A_2$ such that for every even
dimension $2n$ and any 
convex body $K \subset
{\mathbb R}^{2n}$
$$ 
    {\frac {c^Z_{lin}(K)} {c(B^{2n})}}  \leq A_2 \| {\rm Rad} \|^2_{X_K}
 \biggl (
{\frac {{\rm Vol} (K)} {{\rm Vol} (B^{2n})} } \biggr )^{ 1/ n}$$
where $X_K$ is the Banach space whose unit ball is the body $K$.
\end{theorem}

There are wide classes of Banach spaces which are $K$-convex, and
this implies that so are all of their subspaces, with a uniform
bound on their $K$-convexity constant. Thus they will generate
families of convex bodies for which we will have a better bound on
$\gamma_n$ in  Theorem~\ref{cnjV}. In fact, $X$ not being
$K$-convex is equivalent to having an $\varepsilon$-isometric
copy, with respect to the above mentioned Banach-Mazur distance,
of $\ell_1^m$ inside $X$ for every $\varepsilon
> 0$ and every $m$ (see e.g.~\cite{Pi1}). In particular the (infinite dimensional) space
$\ell_p$ of infinite $p$-summable sequences for $1 < p < \infty$ is
$K$-convex (with a $K$-convexity constant depending on $p$) and
therefore for each $1<p<\infty$, for all finite dimensional
$\ell_p^n$ we will have a uniform bound on $\gamma_n$ which will
depend only on $p$. In Section~\ref{IMP} we shall discuss these
examples together with other families of convex bodies for which
Theorem~\ref{The_Main_Theorem-after-rad} holds without the
logarithmic factor (and in particular show that for $\ell_p^n$ we
can have $\gamma_n$ independent also of $p$).

\noindent {\bf Symmetric vs. non-symmetric case:} Since affaine
translations in ${\mathbb R}^{2n}$ are symplectic maps, we shall
assume throughout the text that any convex body $K$ has the origin
in its interior. Moreover, it would be enough to assume in what
follows that our body $K$ is centrally symmetric i.e., $K = -K$.
Indeed, assume we have a general convex body $K\subset \R^{2n}$,
we consider the difference body $K' = K-K$, that is, $K'=\{ x -y :
x,y \in K\}$. Of course, $K \subset K'$, so that for any
symplectic capacity $c(K) \le c(K')$, and $K'$ is centrally
symmetric. Moreover, the Rogers-Shephard inequality~\cite{RogShe}
implies that
\[ {\rm Vol}(K') \le 4^{2n} {\rm Vol} (K),\]
Thus, knowing the inequalities in Theorem~\ref{The_Main_Theorem}
and Theorem~\ref{The_Main_Theorem-after-rad} for $K'$ implies the
same inequalities with an extra factor $16$ for a general convex
body $K$.

\noindent {\bf Non-linear methods:} It is worthwhile to mention that
there are non-linear methods of symplectic embedding constructions,
known as ``symplectic folding" and ``symplectic wrapping", which one
might use when approaching Conjecture~\ref{The-conjecture}. We refer
the reader to~\cite{LM},~\cite{Sch} for more details on this
subject. These methods have been successfully used in some
symplectic embedding constructions of concrete convex bodies, see
e.g.~\cite{Sch},~\cite{T}.

\noindent {\bf Notations:} In this paper the letters
$A_0,A_1,\ldots$ are used to denote universal positive constants
which do not depend on the dimension nor the body involved.
In what follows we identify ${\mathbb R}^{2n}$ with ${\mathbb C}^n$
by associating to $z = x + iy$, where $x,y \in {\mathbb R}^n$, the
vector $(x_1,y_1,\ldots,x_n,y_n)$, and consider the standard complex
structure given by complex multiplication by $i$, i.e.,
$i(x_1,y_1,\ldots,x_n,y_n) = (-y_1,x_1,\ldots,-y_n,x_n)$. Note that
under this identification $\omega(v,iv) = \langle v , v \rangle $,
where $\langle \cdot , \cdot \rangle$ is the standard Euclidean
inner product on ${\mathbb R}^{2n}$.
We will use the notion ``holomorphic plane" for a real
$2$-dimensional plane generated by two vectors of the form $v ,
iv$. We shall denote $P_EK$ the orthogonal projection of a body
$K$ on a subspace $E$, and by $\alpha B_E \equiv B_E(\alpha)$ the
open ball of radius $\alpha$ in the subspace $E$. We will denote
the volume of the $2n$-dimensional Euclidean ball by
$\kappa_{2n}$, and sometimes use the estimates \begin{equation}
\label{Volume-of-the-ball-estimate} \left (\frac{\pi e}{n} \right
)^n \Bigl (\frac{1}{\sqrt{2\pi n}e} \Bigr ) \le \kappa_{2n} =
{\frac {\pi^{n}} {n!}} \le \left ({\frac {\pi e} {n}} \right )^n
\Bigl ({\frac 1  {\sqrt{2\pi n}}} \Bigr ) \le \left ({\frac {\pi
e} {n}} \right )^n. \end{equation}
We will use the notation $\ell_p$ with $1\le p <\infty$
for the space of infinite $p$-summable sequences endowed with the
norm
$(\sum_{i} |x_i|^p)^{1/p}$,  and $\ell_{\infty}$ for
the space of bounded infinite sequences with norm given by
$\sup_{i} |x_i|$. We use $\ell_p^n$ to denote $\R^n$ with
the norm $(\sum_{i=1}^n |x_i|^p)^{1/p}$ for $1\le p<\infty$ and
$\sup_{i=1, \ldots, n} |x_i|$ for $p=\infty$.

\noindent {\bf Structure of the paper:} The paper is organized as
follows. We first consider some concrete examples, where we can
show directly the validity of Conjecture~\ref{The-conjecture}
up to a universal constant. In Section~\ref{Main-ing} we introduce
some of the main ingredients in the proof of
Theorem~\ref{The_Main_Theorem-after-rad}.
In Section~\ref{Proof-of-main-result-section} we show how the method
that works in the concrete examples of Section \ref{ce} enables us
to use some simple linear algebra, together with the bounds
following from the deep work of Pisier mentioned above and the
ingredients from Section~\ref{Main-ing}, to prove
Theorem~\ref{The_Main_Theorem-after-rad}. Finally, in the last
section we discuss some families of convex bodies where
Theorem~\ref{The_Main_Theorem} holds in improved form, i.e., without
the logarithmic factor.

\noindent {\bf Acknowledgment:} We thank Leonid Polterovich for
his support and patience, and for illuminating comments regarding
the text. We thank Felix Schlenk for his helpful advice. The
second named author thanks Alexandru Oancea for sending him his
thesis.

\section{Some Concrete Examples}\label{ce}

In this section we estimate the cylindrical capacity 
of some concrete convex bodies. We show that for these examples
Conjecture~\ref{The-conjecture} holds up to a constant factor.
These examples, in particular the techniques used
for the examples of the ellipsoid and of the cross polytope, will
guide us in the proof of the general case. Other, more general
examples, will be discussed in Section~\ref{IMP}.

\subsection{Ellipsoids}
The fact that for an ellipsoid $E$ one has $c_B(E) = c^Z(E)$ (and
therefore Conjecture~\ref{The-conjecture} holds) 
is well known (see, e.g., \cite{HZ}, \cite{McS}). However, we
wish to remark on its proof since it provides some preliminary
intuition for the general case.

First of all, consider the case where $E$ is a {\em symplectic
ellipsoid}, given by $E = DB^{2n}$ where $D = {\rm
diag}(r_1,r_1,r_2,r_2 \ldots,r_n,r_n)$ with $0<r_1 <r_2<\ldots
<r_n$.
In this case we clearly have that $E$ lies ``between'' a ball and a
cylinder of the same radii, i.e., $ B^{2n}(r_1)\subset E \subset
Z^{2n}(r_1),$ and hence, $c_B(E) = c^Z(E)$.

For $E$ an ellipsoid in general position, one can find suitable
(linear) symplectic coordinates in which $E$ becomes a symplectic
ellipsoid. We discuss this well known fact in Section
\ref{The_Main_Theorem-after-rad}, and as in the proof of Corollary
\ref{decomposition of the l-position} below, it means that
one has $SE =
DB^{2n}$ for some symplectic matrix $S$ and a positive diagonal
matrix $D = {\rm diag}(r_1,r_1,r_2,r_2 \ldots,r_n,r_n)$. Using the
above we see that
\[ c_B(E) = c_B(SE) = c_B(DB^{2n}) = c^{Z}(DB^{2n}) = c^Z(SE) = c^Z(E). \]

\subsection{The Cube}\label{cube}
Let $Q = [-1,1]^{2n}$ be the unit cube in ${\mathbb R}^{2n}$. In
this case it is not hard to check that
$$ c^Z(Q) \leq 2 \pi, \ \ {\rm while} \ \ \biggl ( {\frac { {{\rm Vol}(Q)}} {{\rm
Vol}(B^{2n})}} \biggr )^{1/n} = {\frac 4 {\sqrt[n]{\kappa_{2n}}}}
\ge
 {\frac {4n} {\pi e}}.$$ Indeed, the inequality on the left
 follows since $Q \subset Z^{2n}(\sqrt 2)$ and the inequality on
 the right follows from the estimate~$(\ref{Volume-of-the-ball-estimate})$  above.
%
%
%
Thus, in the case of the cube, we see that Theorem
\ref{cnjV} holds with a constant $\gamma_n \simeq 1/n$, in
particular (at least asymptotically, but it is not hard to check
the constants in general) with constant $1$.
\begin{remark} {\rm
In fact, using that $Q \subset [-1,1]^2\times \R^{2n-2}$, we see
that $c^{Z}(Q) = 4$; however, to embed the cube into $Z^{2n}(r)$
with $\pi r^2 = 4$ one steps out of the linear
category.}\end{remark}

For the following linear image of a cube: ${\widetilde Q} = DQ$,
where $D = diag(a_1,b_1, \ldots,a_n,b_n)$ for some positive real
numbers $a_i,b_i$ we have
\[ c^Z({\widetilde Q}) \leq 2 \pi \min_i a_ib_i \ \ {\rm while}
\ \ 
\biggl ( {\frac {{\rm Vol}({\widetilde Q})} {{\rm Vol}(B^{2n})}}
\biggr )^{1/n} = {\frac {4} {\sqrt[n]{\kappa_{2n}}}
} \sqrt[n]{{ \prod a_ib_i}} \ge {\frac {4n} {\pi e}} \sqrt[n]{{
\prod a_ib_i}}
\]
Indeed,  
the linear symplectic transformation $S = diag \Bigl ( \sqrt
{\frac {b_1} {a_1}},\sqrt {\frac {a_1} {b_1}}, \ldots,\sqrt {\frac
{b_n} {a_n}},\sqrt  {\frac {a_n} {b_n}} \Bigr ) $ satisfies that
$S{\widetilde Q} \subset Z^{2n}(\min {\sqrt { 2 a_i b_i}})$
(Again, the bound can be improved by using a non-linear
symplectomorphism to $c^Z({\widetilde Q}) \le  \min_i 4a_ib_i$).
We see thus that applying a diagonal transformation to the cube
 improves the inequality, and Theorem \ref{cnjV} holds for
this body with a constant $\gamma_n \simeq 1/n$.

\subsection{The Cross-Polytope}\label{crop}

The cross-polytope in dimension $2n$ is the 
polytope corresponding to the convex hull of the $4n$ points
$\pm e_i, i = 1, \ldots, 2n$
formed by permuting the coordinates $(\pm1, 0, 0, ..., 0)$. Since
this is the unit ball of the $\ell_1$ norm on ${\mathbb R}^{2n}$,
we shall denote it by $B^{2n}_1$, so,
$B^{2n}_1 = {\rm Conv}\{ \pm e_i\}$.
We claim that
\[ c^Z(B^{2n}_1) \le {\frac {\pi} n} \ \ {\rm and} \ \
  \biggl ( {\frac {{\rm Vol}(B^{2n}_1)} {{\rm
Vol}(B^{2n})}} \biggr )^{1/n} =  \left ( {\frac {{\frac {2^{2n}}
{2n!}}} {{\kappa}_{2n} }} \right )^{1/n}  \ge {\frac 2 {\pi n}}.
\]
%
%
%

The inequality on the right hand side follows from a direct
computation (using the above estimate
$(\ref{Volume-of-the-ball-estimate})$). In order to estimate from
above the cylindrical capacity of $B^{2n}_1$ we will find a
holomorphic $2$-dimensional plane $E$ such that the projection
$P_E B^{2n}_1$ of $B^{2n}_1$ into $E$ is ``small". More precisely,
consider the following unit vectors
\[ v = {\frac 1 {\sqrt {2n}}}(1,\ldots,1), \ \ iv = {\frac 1
{\sqrt {2n}}}(-1,1,\ldots,-1,1).\]
Note that $v \perp iv$ and that $\omega_{std}(v,iv)=1$. Denote $E_0
= {\rm span} \{e_1,e_2 \}$ and $ E = {\rm span} \{ v,iv \}$. Note
that $E = UE_0$ for some $U \in U(n)$. Next we consider the
projection of $B^{2n}_1$ to the subspace $E$. It follows from the
definition that
\[P_E B^{2n}_1 = P_E {\rm Conv} \{ \pm e_i \} = {\rm Conv} P_E
\{\pm e_i\}.$$ A direct computation shows that
$$P_E e_{j} = \langle e_j,v \rangle v + \langle e_j,iv \rangle iv
=
  \begin{cases}
     {\frac 1
n}(1,0,\ldots,1,0) & j~\text{is odd}, \\
     {\frac 1 n}(0,1,\ldots,0,1) & j~\text{is even}.
  \end{cases}
\]
For every $j$ we have that the Euclidean norm $\| P_E e_j \| =
{\frac {1} {\sqrt{n}}}$, and thus the diameter of $P_E B^{2n}_1$ is
equal to ${\frac {1} {\sqrt{n}}}$ which in turn implies that
$c^Z(B^{2n}_1) = c^Z(U^{-1}B^{2n}_1) \leq {\frac \pi n}$.

Moreover, a direct calculation shows that the cross-polytope
$B^{2n}_1$ includes the Euclidean ball ${\frac {1}
{\sqrt{2n}}}B^{2n}$, so that $c_B(B^{2n}_1)\ge {\frac {\pi}
{2n}}$. In particular, we get that up to constant $2$ the two
capacities $c_B$ and $c_Q$ are equivalent, and hence
Theorem~\ref{cnjV} holds with $\gamma = 2$. By using the
bound~$(\ref{Volume-of-the-ball-estimate})$ for the volume of the
Euclidean ball we get the slightly better constant $\gamma = e/\pi$.

Next we consider the linear image  ${\widetilde B^{2n}_1} = D
B^{2n}_1$ where $D = diag(a_1,b_1, \ldots,a_n,b_n)$ for some
positive real numbers $a_i,b_i$ such that $\det D =1$. By applying
the linear symplectomorphism $S$ from the above example of the
cube we can assume without loss of generality that $a_i=b_i$ for
$i=1, \ldots ,n$. As before we are looking for a $2$-dimensional
holomorphic plane $E$ such that the projection of ${\widetilde
B^{2n}_1}$ into $E$ has small diameter. We choose the direction
\[ {\widehat v} =  \gamma ( a_1^{-1},a_1^{-1},\ldots,a_n^{-1},a_n^{-1}), \ \
{\rm where} \ \ \gamma = \Bigl ( \sum_{i=1}^n { 2 {a_i^{-2}}}
\Bigr )^{-{\frac {1} 2}}, \]
which is easily checked to be a direction on which the
projection of ${\widetilde B^{2n}_1}$ has minimal length. Together
with it we take
$ i{\widehat v} =  \gamma (
-a_1^{-1},a_1^{-1},\ldots,-a_n^{-1},a_n^{-1} )$ (both are unit
vectors). Note that ${\widehat v} \perp i{\widehat v}$ and
$\omega_{std}({\widehat v},i{\widehat v})=1$. Let $ E = {\rm span}
\{ {\widehat v},i{\widehat v} \}$. In order to bound the diameter
of the projection of ${\widetilde B^{2n}_1}$ into the subspace $E$
we compute the projections:
\[
  \begin{cases}
   P_E a_k e_{2k-1} = a_k \langle e_{2k-1},{\widehat v} \rangle {\widehat v} + a_k \langle
 e_{2k-1},i{\widehat v} \rangle i {\widehat v} = \gamma^2
(2a_1^{-1},0,\ldots,2a_n^{-1},0)  &  \\
P_E a_k e_{2k} = a_k \langle e_{2k}, {\widehat v} \rangle
{\widehat v} + a_k \langle e_{2k}, i {\widehat v} \rangle i
{\widehat v} = \gamma^2 (0,2a_1^{-1},\ldots,0,2a_n^{-1}) &
  \end{cases} \]
It follows that the Euclidean lengths $\|P_E a_k e_{2k-1}\|=\|P_E
a_k e_{2k}\| = \sqrt{2}\gamma$ for $k=1,\ldots ,n$. We notice that
$n(\sqrt{2}\gamma)^2$ is the harmonic mean of the numbers $a_i, i =
1, \ldots, n$, and therefore is smaller than their geometric mean
which equals 1. Thus, $\gamma \sqrt{2}\le 1/\sqrt{n}$, and hence, as
in the not-distorted $B_1^{2n}$, we have that $c^Z({\widetilde
B^{2n}_1}) \leq {\frac \pi n}$. So, in this case, again,
Theorem~\ref{The_Main_Theorem} holds with a constant instead of a
logarithmic factor.

\section{First estimates} \label{Main-ing}

One of the main ingredients in the proof of
Theorem~\ref{The_Main_Theorem-after-rad} is a connection between the
symplectic measure of a convex body $K$, given by its symplectic
capacity, and classical notions of the ``width'' of a convex body
$K$. To be more precise we need several further definitions. For a
non-empty centrally symmetric convex body $K$ in ${\mathbb R}^{2n}$
we denote by $\| \cdot \|_K$ the norm on ${\mathbb R}^{2n}$ induced
by $K$, that is, $\|x \|_K = \inf\{ r: x\in rK\}$. We set
\[ M(K) = \int_{S^{2n-1}} \|x\|_K \sigma(dx),\]
for the average of the norm $\| \cdot \|_K$ on the sphere
$S^{2n-1}$, and define $M^*(K) :=M(K^{\circ})$ where $K^{\circ}$ is
the polar body of $K$ defined by $K^{\circ} := \{ x \in {\mathbb
R}^{2n} \, :  \,  \langle x, y \rangle  \leq 1, \, \forall y \in K
\}$. The number $M^*(K)$ is called half the mean width of $K$ because
\[ M^*(K) = \int_{S^{2n-1}} \sup_{y\in K}  \langle x, y \rangle  \sigma(dx), \]
where we integrate over all unit directions $x$ half the distance
between two parallel hyperplanes touching $K$ and perpendicular to
the vector $x$ (which is called  the width of $K$ in direction $x$).
We may assume without loss of generality that $K$ is indeed
centrally symmetric, as was explained in the introduction.

We will sometimes prefer to use instead of the mean width $M^*(K)$
a discrete analogue called the (normalized) Rademacher average,
which is defined as 
\[ s^*(K) = {\rm Ave}_{\varepsilon_i = \pm {\frac 1 {\sqrt{2n}}}}
\sup \bigl \{ \sum_i \varepsilon_i x_i : x \in K \bigr \} = {\rm
Ave}_{\varepsilon \in \{  {\frac {-1} {\sqrt {2n}}}, {\frac 1
{\sqrt{2n}}} \}^{2n} }\| \varepsilon \|_{K^{\circ}}
\]

\noindent That is, we average the dual norm of the body $K$ (or half
the width of $K$), not on the whole sphere of directions as in the
definition of the mean width $M^*(K)$, but only on the vertices of
the normalized cube (in Banach Space theory one usually uses the
non-normalized version, $r^* = \sqrt{2n} s^*$).

This parameter is much more similar to $M^*$ than it seems at
first. Indeed one can show (see \cite{MS}) that there exists a
universal constant $A_3$ such that for every dimension $2n$ and
every symmetric convex body $K\subset \R^{2n}$ one has
\begin{equation} \label{r*-M*-relation}
A_3^{-1} s^*(K) \leq M^*(K) \leq A_3 \|{\rm
Rad}\|_{X_K} s^*(K),
\end{equation}
where $X_K$ is the Banach space with unit ball $K$. Below we use
only the left hand side inequality, which is not difficult to prove,
and is true with constant $A_3 = \sqrt{\pi/2}$. In fact, we could
instead use below the trivial fact that $\int_{U(n)}s^*(UK) =
M^*(K)$. The following theorems give an upper bound for the
cylindrical capacity of a convex body $K$ in terms of $M^*(K)$ and
$s^*(K)$ respectively. (Notice that by the above remarks
Theorem~\ref{r*estimate} is formally stronger than
Theorem~\ref{M*estimate}. However, they admit an almost identical
proof).


\begin{theorem}\label{M*estimate}
There exists a universal constant $A_4$ such that for every even
dimension $2n$ and any centrally symmetric convex body $K \subset
{\mathbb R}^{2n}$, there exist a holomorphic plane $E \subset
{\mathbb R}^{2n}$ such that
$$ P_E (K) \subset A_4 B^2_E( M^*(K) ).$$
In particular, it follows from the monotonicity property of
symplectic capacities that
\[ c(K) \leq \pi A_4^2 M^*(K)^2.\]
\end{theorem}


\begin{theorem}\label{r*estimate}
There exists a universal constant $A_5$ such that for every even
dimension $2n$ and any centrally symmetric convex body $K \subset
{\mathbb R}^{2n}$, there exist a holomorphic plane $E \subset
{\mathbb R}^{2n}$ such that
$$ P_E (K) \subset A_5 B^2_E( s^*(K) ).$$
In particular, it follows from the monotonicity property of
symplectic capacities that
\[ c(K) \leq \pi A_5^2 s^*(K)^2.\]
\end{theorem}

Before we prove these two theorems let us discuss the relation
between Theorem~\ref{M*estimate} and
Theorem~\ref{The_Main_Theorem-after-rad}. First notice that the
estimate in Theorem~\ref{M*estimate} is weaker than
Conjecture~\ref{The-conjecture}: Indeed, Urysohn's inequality (see
e.g.~Corollary $1.4$ in~\cite{Pi1}) gives
\[ \biggl ({\frac { {{\rm Vol}(K)}} {{\rm Vol}(B^{2n})}} \biggr) ^{1 / {2n}}  \le M^*(K) .\]
However, one can ask whether there are bodies for which there is
equivalence of the two.
In some
sense this is indeed the case. In particular, there exists a
universal constant $A_6$ such that every symmetric
convex body $K\subset {\mathbb R}^{2n}$ has a position $K'$ (that
is, a volume preserving linear transformation $T$ of ${\mathbb
R}^{2n}$ and $K' = TK$) in which
\[ M^*(K') \le
A_6 \| {\rm Rad} \|_{X_K} \biggl ({\frac { {{\rm Vol}(K)}} {{\rm
Vol}(B^{2n})}} \biggr) ^{ 1 / {2n}}, \] where $X_K$ is the Banach
space whose unit ball is $K$. We will discuss this well known fact
in more detail in Section~\ref{Proof-of-main-result-section}.
Recall that for a body $K\subset \R^{2n}$ we have $\| {\rm Rad}
\|_{X_K} \le c\log(2n+1)$ for a universal $c$. Unfortunately, the
above mentioned transformation $T$ need not be symplectic, and we
will address this problem too in
Section~\ref{Proof-of-main-result-section}. Further, for some
classes of bodies, the term $\| {\rm Rad} \|_{X_K}$ above can be
eliminated. We discuss this in Section \ref{IMP}.

In the remainder of this section we
prove Theorems~\ref{M*estimate} and~\ref{r*estimate}.
Recall that Markov's
inequality states that if $(X,S,\mu)$ is a measure space, $f$ is a
measurable real-valued function, and $t > 0$, then
\[ \mu ( \{x \in X \ | \ |f(x)| \geq t \}) \leq {\frac 1 t} \int_X
|f|d\mu.\]

\begin{proof}[{\bf Proof of Theorem~\ref{M*estimate}}]
Consider the unit sphere ${S}^{2n-1} \subset {\mathbb R}^{2n}$
equipped with the canonically defined normalized Haar measure
$\mu$. Define $f : {S}^{2n-1} \rightarrow {\mathbb R}$ by $f(x) =
\|x\|_{K^{\circ}}$. It follows from Markov's inequality above that
for $t = \alpha M^*(K)$
\begin{equation} \label{eq1} \mu ( \{x \in
{S}^{2n-1} \ \Big | \ \| x \|_{K^{\circ}} \geq  \alpha M^*(K) \})
\leq {\frac 1 {\alpha M^*(K)}} \int_{{S}^{2n-1}}  \|x
\|_{K^{\circ}} d\mu = {1/\alpha}. \end{equation}
Next, for every unit vector $x =
(x_1,x_2,x_3,x_4,\ldots,x_{2n-1},x_{2n}) \in {S}^{2n-1}$
consider $i x = (-x_2,x_1,-x_4,x_3,\ldots,-x_{2n},x_{2n-1}) \in
{S}^{2n-1}$.
 Note that $x \perp ix$ and
$\omega(x,ix) = 1$. By substituting $\alpha$ large enough, say $3$,
in~($\ref{eq1}$) above, we get that for at least measure $1-1/\alpha
= 2/3$ of the vectors $x$ on the sphere we have that $\|x
\|_{K^{\circ}} \leq \alpha M^*(K)$. So, at least $1-2/\alpha = 1/3$
of the couples $(x,ix)$ satisfy that this is true for {\em both} of
them. Hence, for at least $1/3$ of the vectors $x \in {S}^{2n-1}$ we
have
\[ P_EK \subset \alpha {\sqrt 2 }  M^*(K)B_E^2,
\ \ {\rm where} \ E = {\rm span} \{x,ix \}.\]
%
Note that the subspace $E$ is a unitary image
of the subspace $E_0 = {\rm span} \{e_1,e_2\},$ where $e_1 =
(1,0,\ldots,0)$, $e_2 = (0,1,\ldots,0)$. Since a unitary
transformation preserves the symplectic structure and since the
symplectic capacity is monotone, we have that the capacity
$c^Z(K)$ is at most the capacity of a ball of radius $ \alpha
{\sqrt 2}
 M^*(K)$. This completes the proof of the theorem, with $A_4 = 3\sqrt{2}$.
\end{proof}

\begin{remark}
{\rm Actually, this theorem can be deduced directly from Dvoretzky's
theorem (see e.g.~\cite{MS}) about random sections or projections of
convex bodies which implies, roughly speaking, that most projections
of a convex body $K$  of the appropriate dimension or lower are
Euclidean balls of diameter approximately $M^*(K)$. This can be made
precise: for $2$-dimensional sections given by $UE_0$ for $U\in
U(n)$ and $E_0 = span \{ e_1, e_2\}$, it can be shown that under
mild assumptions the above is true for large (Haar) measure of $U\in
U(n)$ (see, e.g., \cite{AA}). However, since all we need is a
$2$-dimensional subspace, and since we are willing to sacrifice
universal constants (such as $\sqrt{2}$ above), a much more low-tech
approach based only on Markov's inequality worked equally well.}
\end{remark}

\begin{proof}[{\bf Proof of Theorem~\ref{r*estimate}}]
This proof is almost identical to the above. Instead of considering the unit
sphere we consider $\{-1/\sqrt{2n}, 1/\sqrt{2n}\}^{2n} \subset
{\mathbb R}^{2n}$ equipped with the uniform (normalized) counting
measure $\mu$. It follows from Markov's inequality that for $t =
\alpha s^*(K)$
\begin{equation} \label{eq1r} \mu ( \{\eps \in
\{\pm 1/\sqrt{2n} \}^{2n} \ \Big | \ \| \eps \|_{K^{\circ}} \geq
\alpha s^*(K) \})  \le {1/\alpha}. \end{equation}
As before, for every  vector $x =
(x_1,x_2,x_3,x_4,\ldots,x_{2n-1},x_{2n}) \in
\{\pm1/\sqrt{2n}\}^{2n}$ we consider $ix =
(-x_2,x_1,-x_4,x_3,\ldots,-x_{2n},x_{2n-1})$, which is also in
$\{\pm1/\sqrt{2n}\}^{2n}$. Again $x \perp ix$ and $\omega(x,ix) =
1$, and with $\alpha=3$ in~($\ref{eq1r}$) above we get that at least
$1-2/\alpha = 1/3$ of the $2^{2n}$ couples $(x,ix)$ satisfy that
this is true for {\em both} of them. In particular we have one such
couple with
\[ P_EK \subset \alpha {\sqrt 2 }  s^*(K)B_E^2,
\ \ {\rm where} \ E = {\rm span} \{x,ix \}.\]
As before the subspace $E$ is a unitary image of the subspace $E_0
= {\rm span} \{e_1,e_2\},$ so the capacity $c^Z(K)$ is at
most the capacity of a ball of radius $ \alpha {\sqrt 2} s^*(K)$.
This completes the proof of the theorem, with $A_5 = 3\sqrt{2}$.
\end{proof}

\section{
Proof of Theorem~\ref{The_Main_Theorem-after-rad}}
\label{Proof-of-main-result-section}

As mentioned before, a key ingredient in the proof of
Theorem~\ref{The_Main_Theorem-after-rad} is the upper estimate of
the symplectic capacity of a convex body $K$ in terms of its mean
width i.e.,
$$ c(K) \leq \pi (A_4 M^*(K))^2$$
given by Theorem~\ref{M*estimate}. In order to find an upper bound
for the mean width of a convex body in terms of its volume radius
we use a result by Figiel and Tomczak-Jaegermann~\cite{FT}, which
uses a previous result of Lewis~\cite{Le}, stating that every
centrally symmetric convex body $K \subset {\mathbb R}^n$ has a
position $TK$, where $T$ is a volume preserving linear
transformation, such that
$$M(TK)M^*(TK) \leq A_6 \| {\rm Rad} \|_{X_K},$$ where $X_K$ is
the Banach space whose unit ball is the body $K$ and $A_6$ is
universal. Combining this with the fact that
\[ {\frac {1} {M(TK)}} \leq \biggl ({\frac { {{\rm Vol}(TK)}}
{{\rm Vol}(B^{n})}} \biggr) ^{1/n}, \]
where $B^n$ is the $n$-dimensional Euclidean unit ball,  which follows from
polar integration and H\"{o}lder's inequality,
we conclude that
\begin{theorem} \label{l-position theorem-after-rad} There
exists a universal constant $A_6$ such that for every symmetric
convex body $K\subset {\mathbb R}^{2n}$ there exists a position
$TK$, where $T$ is a volume-preserving linear transformation, for
which
\[ M^*(TK) \leq A_6 \| {\rm Rad} \|_{X_K}  \,
\biggl ({\frac { {{\rm Vol}(K)}} {{\rm Vol}(B^{2n})}} \biggr) ^{1/{2n}} \] 
\end{theorem}
This is already close to our goal. However, it is important to
note that the above mentioned transformation $T$ need not be
symplectic, and hence Theorem~\ref{The_Main_Theorem-after-rad}
does not follow directly from the combination of
Theorem~\ref{M*estimate} and Theorem~\ref{l-position
theorem-after-rad}. However, these theorems serve as a tool and as
motivation for the line of proof of
Theorem~\ref{The_Main_Theorem-after-rad}. Our next step is
therefore to deal with the ``non-symplectivity"  of the
transformation $T$. To this end, we shall need the following
well known fact about the simultaneously normalization of a
symplectic form and an inner product, which we already used in the
example of an ellipsoid in Section~\ref{ce} (see e.g.~\cite{McS},
page $57$).

\begin{lemma}\label{fp2} Let $(V,\omega)$ be a symplectic vector space and
let $g : V \times V \rightarrow {\mathbb R}$ be an inner product.
Then there exists a basis $\{u_1, v_1,\ldots, u_n,v_n\}$ of
$V$ which is both $g$-orthogonal and $\omega$-standard, that is,
$g(v_i, u_j)= 0$, $g(u_i, u_j) = c_j\delta_{i,j}$,
$g(v_i, v_j) =d_j\delta_{i,j}$, and
$\omega(u_i, u_j) = \omega(v_i, v_j) = 0$, $\omega(u_i, v_j) = \delta_{i,j}$.
Moreover, this basis can be chosen such that $c_j=d_j$ for all
$j$.
\end{lemma}

Straightforward linear algebra gives

\begin{corollary} \label{decomposition of the
l-position} Let ${\mathbb R}^{2n}$ be equipped with the standard
symplectic structure and the standard inner product. Let $T$ be a
volume preserving $2n$-dimensional real matrix. Then there exists
a linear symplectic matrix $S \in {\rm Sp}({\mathbb R}^{2n})$ and
an orthogonal transformation $W \in O(2n)$ such that
\[ T  = WDS, \ \ {\rm where} \ \ \ 0 < D = {\rm
diag}(r_1,r_1,r_2,r_2,
\ldots,r_n,r_n), \ {\rm and} \  \prod_i r_i = 1. \]
\end{corollary}

We postpone the proof of Corollary~\ref{decomposition of the
l-position} to the end of this section. Combining
Equation~($\ref{r*-M*-relation}$) with Theorem~\ref{l-position
theorem-after-rad} and Corollary~\ref{decomposition of the
l-position},
and using the fact that $M^*$ is
invariant under orthogonal transformations we get

\begin{theorem} \label{l-position theorem-second-version}
There exist universal constants $A_3$ and $A_6$ such that for any
dimension $2n$, for every symmetric convex body $K\subset {\mathbb
R}^{2n}$ there exists a positive diagonal transformation $D = {\rm
diag}(r_1,r_1,r_2,r_2, \ldots,r_n,r_n)$ with $\prod r_i = 1 $  and a
symplectic position $K' = SK$ where $S$ is a linear symplectic
transformation such that
\[  s^*(DK') \leq A_3 M^*(DK')   \leq A_3A_6 \| {\rm Rad} \|_{X_K}  \biggl ({\frac { {{\rm Vol}(K)}} {{\rm Vol}(B^{2n})}} \biggr) ^{ 1 / {2n}}.\]
\end{theorem}

In order to complete the proof of
Theorem~\ref{The_Main_Theorem-after-rad} we shall need the
following proposition.

\begin{proposition} \label{Main-argument}
There exists a universal constant $A_7$ such that for any convex
body $K\subset \R^{2n}$ and for every diagonal matrix $D = {\rm
diag}(r_1,r_1,r_2,r_2, \ldots,r_n,r_n)$ with $\prod r_i = 1 $ there
exists a holomorphic plane $E={\rm span}\{v,iv\}$ such that the
orthogonal projection of the body $K$ to the subspace $E$ satisfies
$$ P_E(K) \subset A_7 B_E^2(s^*(DK))$$
\end{proposition}

Postponing the proof of Proposition~\ref{Main-argument} we first use
it to prove Theorem~\ref{The_Main_Theorem-after-rad}.
\begin{proof}[{\bf Proof of Theorem~\ref{The_Main_Theorem-after-rad}}]
Let $K \subset {\mathbb R}^{2n}$ be a symmetric convex body. It
follows from Theorem~\ref{l-position theorem-second-version} that
there exists a symplectic linear image $K' = SK$ of $K$ where $S \in
{\rm Sp}({\mathbb R}^{2n})$ and a positive diagonal matrix $D = {\rm
diag}(r_1,r_1,r_2,r_2, \ldots,r_n,r_n)$ with $\prod r_i = 1 $ such
that
\[  s^*(DK') \leq  A_8  \| {\rm Rad} \|_{X_K}  \biggl ({\frac
{ {{\rm Vol}(K)}} {{\rm Vol}(B^{2n})}} \biggr) ^{ 1 / {2n}}
\] where $A_8= A_3A_6$ is universal. Next we apply
Proposition~\ref{Main-argument} to the body $K'$ and to the above
mentioned diagonal matrix $D$, and conclude that there exists a
holomorphic plane $E$ such that
$$ P_E(K') \subset A_7 B_E^2(s^*(DK'))$$
Since $E$ is holomorphic and since
the symplectic capacity is monotone, we have that for
every symplectic capacity $c$ we have
\[ c(K) = c(K') \le c^Z_{lin}(K') \le \pi
A_7^2 A_8^2 \| {\rm Rad} \|^2_{X_K}  \biggl ({\frac { {{\rm
Vol}(K)}} {{\rm Vol}(B^{2n})}} \biggr) ^{ 1 / {n}},\] and the
proof of Theorem \ref{The_Main_Theorem-after-rad} is complete.
\end{proof}

The rest of this section is devoted to the proofs of
Proposition~\ref{Main-argument} and of Corollary~\ref{decomposition
of the l-position}.

\begin{proof}[{\bf Proof of Proposition~\ref{Main-argument}}]
Let $K$ be a centrally symmetric convex body and let $D = {\rm
diag}(r_1,r_1,\ldots,r_n,r_n)$ be a positive diagonal matrix with $
\det D = 1$. It follows from the proof of Theorem~\ref{r*estimate}
applied to the body $DK$ that there exists an holomorphic plane
${\widehat E} = {\rm span} \{v,iv \}$, where $v$ and $iv$ are
vertices of the cube $\{ {\frac 1 {\sqrt{2n}}}, {\frac {-1}
{\sqrt{2n}}} \}^{2n} $ such that
\begin{equation} \label{conclusion-from-r^*}
\|v\|_{(DK)^{\circ}} = \sup_{z \in DK} | \langle v , z \rangle| \leq
A_5 s^*(DK), \  \  \|iv\|_{(DK)^{\circ}} = \sup_{z \in DK} | \langle iv
, z \rangle| \leq  A_5 s^*(DK) \end{equation}
Next, denote by $\| \cdot \|_2$ the Euclidean norm and consider the
vectors $$v' = {\frac {Dv} {\|Dv \|_2}} = \Bigl
({n\over\sum_{i=1}^nr_i^2} \Bigr )^{1/2}
(v_1r_1,v_2r_1,v_3r_2,v_4r_2,\ldots,v_{2n-1}r_n,v_{2n}r_n)$$ and
$$ i v' = {\frac {iDv} { \| D
v \|_2 }}    = {\frac {D i v} { \| D  v \|_2 }} = \Bigl({n\over
\sum_{i=1}^nr_i^2}\Bigr)^{1/2}
(-v_2r_1,v_1r_1,-v_4r_2,v_3r_2,\ldots,-v_{2n}r_n,v_{2n-1}r_n)$$ Note
that these two unit vectors satisfy of course that $v' \perp iv'$
and $\omega_{std}(v',iv') = 1$. Moreover by the geometric arithmetic
mean inequality we see that
\begin{equation}
\label{GAmean} \|Dv\|_2 = \Bigl ({1\over n} \sum_{i=1}^nr_i^2
\Bigr )^{1/2} \geq \Bigl ( \prod_{i=1}^n r_i^2 \Bigr )^{ 1 / {2n}}
= 1 . \end{equation} Next, denote $E = {\rm span} \{v',iv'\}$. A
straightforward computation shows that
\begin{eqnarray*} P_{E}K  & = & \{ \alpha v' + \beta iv' \ | \ \alpha
= \langle v',y \rangle , \ \beta = \langle i v',y\rangle , \ \
{\rm where} \ y \in K \} \\
& = & \{ \alpha v' + \beta iv' \ | \ \alpha = {\frac 1 {\|Dv\|_2}
}\langle v, Dy \rangle , \ \beta = {\frac 1 {\|Dv\|_2}} \langle
  iv,Dy \rangle , \ \ {\rm where} \ y \in K \} \\ & = & \{ \alpha v' + \beta iv' \ | \ \alpha
= {\frac 1 {\|Dv\|_2} }\langle v, z \rangle , \ \beta = {\frac 1
{\|Dv\|_2}} \langle
  iv,z \rangle , \ \ {\rm where} \ z \in DK \} \\
  & \subseteq   & \{ \alpha v' + \beta iv' \ | \ \alpha
=  \langle v, z \rangle , \ \beta =  \langle
  iv,z \rangle , \ \ {\rm where} \ z \in DK \},
\end{eqnarray*}
where the last inclusion follows from inequality~($\ref{GAmean}$)
above. Combining this with Equation~($\ref{conclusion-from-r^*}$)
above we get that
$$ P_E(K)  \subseteq  {\sqrt 2} A_5 s^*(DK) B_E^2. $$
The proof of the proposition is now complete with $A_7 = {\sqrt
2}A_5$.
\end{proof}

\begin{proof}[{\bf Proof of Corollary~\ref{decomposition of the
l-position}}] It follows from the Lemma \ref{fp2} that there exists
a linear symplectic matrix ${ S}$ and a positive diagonal matrix ${
D^2}={\rm diag}({ r_1^2},{ r_1^2},{ r_2^2},{ r_2^2}, \ldots,{ r_n^2}
,{ r_n^2})$ such that $T^tT = { S}^t{ D^2} { S}$ (since $T^tT$ is
symmetric and positive definite). Note that for the ellipsoid $E =
T^{-1} B^{2n}$ we have
\begin{eqnarray*} E & := &  \{ x \ | \ \langle x,T^tTx \rangle  \leq 1 \} = \{ x \ |
\ \langle x,{ S}^t { D^2} { S}x\rangle  \leq 1 \} = { S}^{-1} \{ y
\ | \ \langle { S}^{-1}y, { S}^t { D^2} y\rangle  \leq 1 \}
\\ & = & { S}^{-1} D^{-1} \{ z \ | \ \langle { S}^{-1}
D^{-1}z, { S}^t Dz\rangle  \leq 1 \}  =  ({D}S)^{-1} \{ z \ | \
\langle S^{-1} D^{-1}z,({D}S)^tz\rangle  \leq 1 \} \\ & = &  (
DS)^{-1} \{ z \ | \  \langle DS ( DS)^{-1}z,z\rangle  \leq 1 \} =
({ D}S)^{-1}B^{2n}
\end{eqnarray*}
Thus we get that $B^{2n} = TE$ and $DS(E) = B^{2n}$, so we conclude
that there exists an orthogonal transformation $W \in O(2n)$ such
that $T({ D}S)^{-1}= W$. We thus conclude that ${T} = W D S$, as
stated.
\end{proof}

\section{Improvements for special families of  bodies}\label{IMP}

In this section we describe wide classes of convex bodies for which
the logarithmic term in Theorem~\ref{The_Main_Theorem} can be
disposed of.

We begin with describing the example of zonoids. The class of
zonoids consists of symmetric convex bodies which can be
approximated, in the Hausdorff sense, by Minkowski sums of line
segments. Bodies which are Minkowski sums of segments are called
zonotopes, and any zonotope in $\R^n$ can be realized as a linear
image of an $m$-dimensional cube for some $m$.  An example for a
zonoid is the Euclidean ball in $\R^n$, which can be approximated in
the Hausdorff distance up to $\eps$ by (random) orthogonal
projections of $m$ dimensional cubes for $m= C(\eps)n$ (this follows
from Dvoretzky's Theorem, see \cite{MS}). In the paper \cite{GMR} it
was shown that every zonoid $Z\subset {\mathbb R}^{2n}$ with volume
$1$ has some linear image $TZ$ for $T \in SL({\mathbb R}^{2n})$ such
that $M^*(TZ)\le A_9 \sqrt{n}$ for some universal constant $A_9$.
After normalization we see that this means $M^*(TZ)\le A_9' \bigl
({\frac { {{\rm Vol}(Z)}} {{\rm Vol}(B^{2n})}} \bigr) ^{1/{2n}} $.
Thus, applying the same methods as in the proof of
Theorem~\ref{The_Main_Theorem} we conclude that for zonoids the
estimate in Theorem~\ref{The_Main_Theorem} is true without the
logarithmic factor.

Next, we consider the unit ball of $\ell_p^{2n}$ for $1\le p\le
\infty$, denoted $B_p^{2n}$. The case of $p = 1, \infty$ was
discussed in Section~\ref{ce} above, where it was shown that
Theorem~\ref{The_Main_Theorem} holds with a universal constant
instead of the logarithmic factor. For $1 < p < \infty $, one
option, which we {\em do not} use, is to use the well known
computations for $M^*(B_p^{2n})$ (see e.g. \cite{MS}) and for
${\rm Vol}(B_p^{2n})$,
\[ {\rm  Vol} (B_p^{2n}) = {(2\Gamma ({1\over p} +1))^{2n}\over
\Gamma({2n\over p} +1)}.   \]
In what follows we take a different, more geometric approach, and
avoid using the above mentioned estimates.

We begin with the case of $1<p \le 2$, where we can invoke the
following inclusion which is easy to check
\[ (2n)^{{{1/2} -{1/ p}}}B^{2n} \subset B_p^{2n} \subset {(2n)}^{1-1/p}B_1^{2n}. \]
The left hand side inclusion implies that
\[ \biggl ( {\frac {{\rm
Vol}(B_p^{2n})} {{\rm Vol}(B^{2n})}} \biggr )^{1/{n}} \ge
(2n)^{1-2/p} ,\] and the right have side inclusion together with
the results for $B_1^{2n}$ from Subsection~\ref{crop} implies that
 \[c^Z(B_p^{2n})\le (2n)^{2-2/p}c^Z(B_1^{2n})\le
2\pi (2n)^{1-2/p}. \] This completes the case $1<p \le 2$.
Moreover, we see that the constant $2$ we get is universal and
does not depend on $p$.

We turn to the case $p>2$. Here we are even better off, because
we will use the estimates for the cube from Subsection~\ref{cube}
where there was a big difference between the capacity and the
ratio of volumes. Indeed, we may use the inclusions
for $2<p<\infty$ (where the cube is denoted now by $B_{\infty}^{2n}$)
\[ B^{2n} \subset B_p^{2n} \subset B_{\infty}^{2n}. \]
Therefore
\[ \biggl ( {\frac {{\rm
Vol}(B_p^{2n})} {{\rm Vol}(B^{2n})}} \biggr )^{1/{n}} \ge 1 ,\]
and using the results of  Subsection~\ref{cube} we also have
\[ c^Z(B_p^{2n}) \le c^Z(B_{\infty}^{2n}) \le 4, \]
which completes the case $2<p<\infty$ with constant ${\frac 4
{\pi}}$. Had we used the exact estimates for $M^*$ and for the
volume we would have gotten that the theorem holds with a much
better constant (equal to $1$, and getting smaller as $p$ grows).
We remark again that we arrived at a universal constant
independent of both $n$ and $p$.

Other examples for which the logarithmic factor can be omitted
arise from Theorem \ref{The_Main_Theorem-after-rad}. Indeed,
whenever we have good bounds on $\|{\rm Rad}\|_X$ we can replace
the logarithmic factor by these bounds. We could have applied this
scheme in the $\ell_p^{2n}$ case, however that would {\em not}
guarantee that the constants do not depend on $p$, and in fact,
since for $X = \ell_1, \ell_{\infty}$ we know that $\|{\rm
Rad}\|_X$ is not bounded, we would have gotten constants that
depend on $p$ and explode as $p\to 1, \infty$.

We emphasize that $K$-convexity is an infinite dimensional notion,
and indeed since we are concerned with the dependence on dimension
in Theorem \ref{The_Main_Theorem}, if we wish to prove a stronger
bound in special cases then these particular cases have to be
families of convex bodies in dimension tending to infinity. Some
such examples are concrete bodies such as $\ell_p^{2n}$, but another
way to construct such families is to look at some infinite
dimensional Banach space and to consider for example all of its
finite dimensional subspaces and the convex bodies which are their
unit balls. In this case $K$-convexity of the original space $X$
will promise a uniform bound, independent of dimension, in Theorem
\ref{The_Main_Theorem}, instead of the logarithmic factor.

For example, the Rademacher constant $\|{\rm Rad}\|_X$ is bounded
for $X = L_p(Y)$ for any (fixed) Banach space $Y$ and
$1<p<\infty$, which is defined for a probability measure $\mu$
on a set $\Omega$
as
\[ L_p(\Omega,Y) = \Bigl \{ f : \Omega \rightarrow Y \ : \
\|f\|_{L_p(\Omega,Y)} = \Bigl ({\int_{\Omega}}\|f\|^p_Y \, d\mu
\Bigr )^{1/p} < \infty \Bigr \}. \]
Thus, all convex bodies which are unit
balls of (finite dimensional) subspaces of the space $L_p(Y)$
will satisfy
the main theorem without the logarithmic factor (but the bound
will depend on the $K$-convexity constant of $L_p(Y)$).

More generally, $X$ is $K$-convex if and only if $X$ is of
type-$p$ for some $p>1$. The definition of type-$p$ is as follows
(see, e.g. \cite{Pi1}): A Banach space $X$ is called of type-$p$
for $1\le p \le 2$ if there is a constant $C$ such that for all
$m$ and all $x_1, \ldots, x_m \in X$ we have
\begin{equation}\label{tyyppe}
\| \sum_{i=1}^m r_i x_i \|_{L_2(X)} \le C \left(\sum_{i=1}^m
\|x_i\|^p \right)^{1/p}. \end{equation} The smallest constant $C$
for which this holds is called the type-$p$ constant of $X$, and is
denoted $T_p(X)$. For a Hilbert space obviously $T_2(X) = 1$, and we
remark that Kahane's inequality (see e.g.~\cite{MS}) states that for
$1\le p <\infty$ there are constants $K_p$ depending only on $p$ so
that for any $X$ and any $x_1, \ldots, x_m \in X$
\[ \| \sum r_i x_i \|_{L_1(X)} \le \| \sum r_i x_i \|_{L_p(X)}
\le K_p \| \sum r_i x_i \|_{L_1(X)}, \] so that the notion of
type-$p$ does not depend on the choice of $L_2$-average on the left
hand side of~$(\ref{tyyppe})$. It is easily seen from the triangle
inequality that every Banach space has type 1 with $T_1(X) = 1$, and
it follows easily, say from Khinchine inequality, that no non-zero
normed space has type $p>2$.

The theorem we stated above, that $X$ is $K$-convex if and only if
it has some non-trivial type (i.e., type-$p$ for $p>1$) is due to
Milman and Pisier, see e.g.~Theorem~11.3 in \cite{Pi1}. Many more
equivalent formulations are possible, for example the same theorem
also states that $X$ is $K$-convex if and only if it is locally
$\pi$-Euclidean (for the definitions see \cite{Pi1}).

One more example for convex bodies where this method shows that
the logarithmic factor is not necessary are unit balls of the
Schatten classes $C_p^n$ for $1<p<\infty$. To define the Schatten
class spaces let $u$ be an $n\times n$ matrix, so $u^*u$ is
positive definite and symmetric, thus it is orthogonally
diagonalizable with nonnegative eigenvalues $\lambda_1, \ldots,
\lambda_n$. The Schatten class $C_p^n$ is the $n^2$-dimensional
space of all $n\times n$ real matrices equipped with
the norm
\[ \|u\|_{C_p^n} = \left( \sum_{i=1}^n \lambda_i^{p/2}\right)^{1/p}. \]
(Of course we will consider $C_p^{2n}$ when we want to discuss the
symplectic capacity of the unit ball of this space.)
Tomczak-Jaegermann showed in \cite{TJ} that Schatten classes $C_p^n$
have the same type/cotype properties as $L_p$ spaces, so that we
have a uniform estimate in Theorem \ref{The_Main_Theorem} also for
the unit balls of $C_p^n$ when $1<p<\infty$ (but the constant which
replaces the logarithm can depend on $p$).



\begin{remark} {\rm In a paper by Giannopoulos, Milman and Rudelson \cite{GMR}
they prove a theorem which gives a bound on the minimal $M^*(TK)$
over $T\in SL(n)$, which is slightly better than the simple bound
we used above. Recall that we argued that
\[ \min_T M^*(TK) \le
\min_T (M(TK) M^*(TK)) \left(\frac{{\rm Vol}(K)}{{\rm Vol}(B^{2n})}\right)^{1/2n},\]
and then we used the estimate
$\min_T(M(TK) M^*(TK)) \le A_6 \| {\rm Rad} \|_{X_K}$ for a universal
$A_6$.
They show that
(stating their Theorem~B in dual form): For every body $K$ there
is a position $K'= TK$ for $T\in SL(n)$ such that (for a universal
constant $c$), denoting by $d$ half the diameter of $K'$ (so that
$K' \subset d B^{2n}$) one has
\[ \frac{c{M^*(K')}}{\log
(d/M^*(K'))}  \le \left(\frac{{\rm Vol}(K')}{{\rm Vol}(B^{2n})}\right)^{1/2n}.\]
It is easily seen that one
{\em always} has $d(K')/M^*(K') \le C'\sqrt{n}$ (since $K'$
includes a segment of length $2d$, and so $M^*(K') \ge
M^*([-d, d]) = c'/\sqrt{n}$). Therefore their result clearly
implies that there is a position for which
\[ {M^*(K')} \le C_1 (\log C_2 n)
\left(\frac{{\rm Vol}(K')}{{\rm Vol}(B^{2n})}\right)^{1/2n},\] which
is exactly the estimate we derived and used above. However, it is
plausible that in many cases the position $TK$ they use has a better
ratio of diameter and mean width, and then the above estimate give
an improved result. }
\end{remark}

\bigskip
\noindent Shiri Artstein-Avidan  \\
 Department of Mathematics \\
 Princeton University Fine Hall \\
Washington Road, Princeton NJ 08544-1000 USA \\
{\it e-mail}: artstein@math.princeton.edu \\

\bigskip
\noindent
Yaron Ostrover\\
School of Mathematical Sciences \\
Tel Aviv University \\
Tel Aviv 69978, Israel  \\
{\it e-mail}: yaronost@post.tau.ac.il\\


\begin{thebibliography}{}

\bibitem{AA} Artstein-Avidan, S. {\it
A Bernstein-Chernoff deviation inequality and geometric properties of
random families of operators.} To Appear in the Israel Journal (2006).

\bibitem{CHLS} Cieliebak, K., Hofer, H., Latschev, J., Schlenk
F. {\it Quantitative symplectic geometry.} math.SG/0506191.

\bibitem{EH} Ekeland, I., Hofer, H. {\it Symplectic topology
and Hamiltonian dynamics.} Math. Z. {\bf 200} (1989), no. 3,
355--378.

\bibitem{FT} Figiel, T., Tomczak-Jaegermann, N. {\it Projections
onto Hilbertian subspaces of Banach spaces.} Israel J. Math. 33
(1979), 155--171.

\bibitem{GMR} Giannopoulos, A. A., Milman, V.D., Rudelson, M. {\it Convex bodies with minimal mean
width.} Geometric aspects of functional analysis, 81--93, Lecture
Notes in Math., 1745, Springer, Berlin, 2000.

\bibitem{G} Gromov, M. {\it Pseudoholomorphic curves in symplectic manifolds.} Invent. Math. {\bf 82} (1985), no. 2, 307-347.



\bibitem{H} Hermann, D. {\it Non-equivalence of symplectic
capacities for open sets with restricted contact type boundary.}
Pr\'epublication d'Orsay num\'ero 32 (29/4/1998).



\bibitem{Ho} Hofer, H. {\it Symplectic capacities.} Geometry of
low-dimensional manifolds, 2 (Durham, 1989), 15-34, London Math.
Soc. Lecture Note Ser., {\bf 151}, Cambridge Univ. Press,
Cambridge, 1990.

\bibitem{HZ} Hofer, H. and Zehnder, E. {\it Symplectic Invariants
and Hamiltonian Dynamics.} Birkh\"auser, Basel (1994).

\bibitem{J} John, F. {\it Extremum problems with inequalities as
subsidiary conditions.} Courant Anniversary Volume, Interscience,
New York,  187-204 (1948).

\bibitem{L} Lalonde, F. {\it Energy and capacities in symplectic
topology.} Geometric topology (Athens, GA, 1993), 328-374, AMS/IP
Stud. Adv. Math., {\bf 2.1}, Amer. Math. Soc., Providence, RI,
1997.

\bibitem{LM} Lalonde, F., McDuff, D. {\it The geometry of symplectic
energy.} Ann. of Math. 141 (1995) 349-371.

\bibitem{Le} Lewis, D.R. {\it Ellipsoids defined by Banach ideal
norms.} Mathematika, 26:18-29 (1979).

\bibitem{Mc} McDuff, D. {\it Symplectic topology and capacities.}
Prospects in mathematics (Princeton, NJ, 1996), 69-81, Amer. Math.
Soc., Providence, RI, 1999.

\bibitem{McS} McDuff, D. and Salamon, D. {\it Introduction to
Symplectic Topology}, 2nd edition, Oxford University Press, Oxford,
England (1998).

\bibitem{MS} Milman, V.D., Schechtman, G. {\it Asymptotic theory
of finite dimensional normed spaces.} Lectures Notes in Math. 1200,
Springer, Berlin (1986).


\bibitem{Pi1} Pisier, G. {\it The volume of convex bodies and Banach space
geometry.} Cambridge University Press, Cambridge, (1989).


\bibitem{Pi2} Pisier, G. {\it Holomorphic semi-groups and the
geometry of Banach spaces.} Annals of Math. {\bf 115} (1982),
375-392.

\bibitem{RogShe} Rogers, C.A. and Shephard, C. {\it The difference
body of a convex body.} Arch. Math. 8 (1957), 220-233.

\bibitem{Sch} Schlenk, F. {\it Embedding problems in symplectic geometry.} de Gruyter Expositions
in Mathematics, 40. Walter de Gruyter GmbH \& Co. KG, Berlin,
2005.



\bibitem{TJ} Tomczak-Jaegermann, N. {\it The moduli of smoothness and convexity and the
Rademacher averages of trace classes $S_p$ ($1< p < \infty$)},
Studia Math. 50 (1974), 163--182.

\bibitem{T} Traynor, L. {\it Symplectic packing constructions.}
J.Differential Goem. 42 (1995), 411-429.

\bibitem{V} Viterbo, C.
{\it Metric and isoperimetric problems in symplectic geometry.} J.
Amer. Math. Soc. 13 (2000), no. 2, 411--431.

\bibitem{V1} Viterbo, C. {\it Capacit\'es symplectiques
et applications (d'apr\`es Ekeland-Hofer, Gromov).} S\'eminaire
Bourbaki, Vol. 1988/89. Ast\'erisque no. {\bf 177-178} (1989),
Exp. no. 714, 345-362.




\end{thebibliography}
\end{document}